\documentclass[12pt]{amsart}
\usepackage{amsmath, amsthm, amsfonts, amscd, amssymb, latexsym}
\usepackage{color}
\usepackage{enumitem}
\setlist[enumerate,1]{%
label={\normalfont(\arabic*)},ref=\arabic*%
}
\theoremstyle{plain}
\newtheorem{thm}{Theorem}[section]

\theoremstyle{definition}
\newtheorem{defn}[thm]{Definition}
\newtheorem{ex}[thm]{Example}

\numberwithin{equation}{section}
\begin{document}
\title{The second variational formulas for mappings between statistical manifolds}
\title[The second variational formula for statistical manifolds]{ 
The second variational formula of the variational problems for mappings 
between statistical manifolds}
\author{
Hajime Urakawa}
\address{
Graduate School of Information Sciences, 
Division of Mathematics, Tohoku University, 
Aoba 6-3-09, Sendai 980-8579, Japan}
\email{hajime.urakawa.c6@@tohoku.ac.jp}
    \keywords{harmonic map, variational problem, statistical manifold, 
    second variational formula}
  \subjclass[2020]{primary 58E10, secondary 53C42}
  \thanks{
}
\maketitle
\begin{abstract} Recently, H. Furuhata and R. Ueno \cite{1} obtained the first variational formula 
of smooth mappings of a compact statistical manifold to another statistical manifold. 
In this paper, we show the second variational formula of harmonic mappings of a statistical manifold into another one. Furthermore, we define the stability, the index and the nullity for the mappings, and show the weakly stability for the harmonic mappings into any statistical manifolds of non-positive curvature.   
\end{abstract}
\vskip0.6cm\par
\section{Introduction.} 
Let $u\,;M\rightarrow N$ be a smooth map from 
a compact statistical manifold $(M,g,\nabla^M)$ into 
another statistical manifold $(N,h,\nabla^N)$. 
Then, recently, 
H. Furuhata and R. Ueno \cite{1} gave the first variational formula: 
For an arbitrary smooth variation $F=\{u_t\}_{-\epsilon<t<\epsilon}$ of $u$, 
$$
\frac{d}{dt}\Big\vert _{t=0}\,E_2(u_t)=\int_M\langle V,\tau(u)\rangle\,d\mu_g,
$$ 
where 
$$E_2(u)=\frac12\int_M\Vert\tau(u)\Vert^2\,d\mu_g,$$
$d\mu_g$ is the volume element 
of the Riemannian metric $g$, and 
the variation vector field $V$ is given by 
$$
V(x)=\frac{d}{dt}\Big\vert _{t=0}u_t(x), \qquad x\in M. 
$$
\par 
In this paper, we want to give the second variational formula, namely, differentiating 
$E(u)$ twice  i.e., 
$$
\frac{\partial^2}{\partial s\partial t}\Bigg\vert_{(s,t)=(0,0)}E_2(u_{s,t}).
$$
As its application, we derive the Jacobi operator which introducing 
several notions of the Morse theory, 
 the index, the nullity, and the stability of the harmonic mapping. 
 And we give several examples of stable harmonic mappings, and unstable harmonic mappings between statistical manifolds.  
   \vskip0.6cm\par
\section{Preliminaries.}
In this section, we prepare to state the results due to 
H. Furuhata and R. Ueno, the first variational formula 
of the statistic bi-energy (cf. \cite{FU}). 
 Let us consider a smooth map $u:\,(M,g)\longrightarrow (N,h)$ 
 from a compact statistical manifold $(M,g,\nabla^M)$ 
 into another statistical manifold $(N,h,\nabla^N)$. 
 Let $\tau(u)\in \Gamma(u^{-1}TN)$ be the statistical bi-tension field $(1.3)$ 
 for $u$ given by Furuhata and Ueno 
 (cf. \cite{FU}, Pages 2, 10). 
 Namely, they showed the first variation formula as follows. 
 First they defined the functional $E_2(u)$ for every smooth map $u:\,M\rightarrow N$ by 
 \begin{align}
 E_2(u)&:=\frac12\int_M\vert \tau_2(u)\vert{}^2\,d\mu_g, \\
 \frac{d}{dt}\Big\vert _{t=0}E_2(t)&=\int_M\langle V,\tau_2(u)\rangle\,d\mu_g, 
 \end{align}
 where 
 \begin{align}
 \tau_2(u):={\rm tr}_g \Big\{ (X,Y)\mapsto \nabla {}^N_X u_{\ast}Y
 -u_{\ast}{{\nabla}^M_XY} \Big\}
 \in \Gamma(u^{-1}TN).
 \end{align}
 The Euler-Lagrange equation is $\tau_2(u)=0$ where $\tau_2(u)$ is given by 
 
 $$
 \overline{\Delta}^u\tau_2(u)+{\rm div}^g({\rm tr}_gK^M)\tau_2(u)
 -\sum_{i=1}^mL^N(u_{\ast}e_i,\tau_2(u))u_{\ast}e_i-K^N(\tau(u),\tau_2(u)).
 $$
 \par
 In this paper, we see the second variational formula of $E_2(u)$, i.e., 
we calculate explicitly 
 \begin{align}  
 \frac{\partial^2}{\partial s\partial t}\Big\vert _{(s,t)=(0,0)}E_2(u(s,t)), 
 \end{align}
 indeed, for a smooth variation $u_{s,t}$ of 
 with two parameters $s$ and $t$ of a harmonic map $u:\,(M,g,\nabla^M)$ into $(N,h,\nabla^N)$, 
 we determine 
 the following Hessian of the energy $E_2$ at $u$: 
  \begin{align}  
 H(E_2)_u(V,W)=\frac{\partial^2}{\partial s\partial t}\Big\vert _{(s,t)=(0,0)}E_2(u(s,t))
 \end{align}
 which introduces the second order elliptic differential operator $J_u$ 
 acting on the space  $\Gamma(u^{-1}TN)$, called the {\it Jacobi operator} satisfying that 
 \begin{align}
 H(E_2)_u(V,W)=\int_Mh(J_u(V),W)\,d\mu_g,\qquad V,\,\,W\in \Gamma(u^{-1}TN). 
 \end{align}
\section{The second variational formula of the statistical bi-energy.}
In this section, we will state the second variational formula of the statistical bi-energy. 
\par
Let  $u:\,\, (M,g,{\nabla}^M)\rightarrow (N,h,{\nabla}^N)$ be a statistical harmonic map, 
namely, for a variation $u_{s,t}$  with two parameters  $(-\epsilon<s,t<\epsilon)$  
of the statistical harmonic mapping  $u$, 
let 
$F:\, (-\epsilon,\epsilon)\times (-\epsilon,\epsilon)\times M\rightarrow N$ be 
a smooth variation with two parameters
$$
(-\epsilon,\epsilon)\times (-\epsilon,\epsilon)\times M\ni 
(s,t,x)\mapsto F(s,t,x)=u_{s,t}(x)\in N
$$
with 
$F(0,0,x)=u_{0,0}(x)=u(x)\in N$, $x\in M$. 
\par
We want to give the second variation formula the energy $E$ for two vector 
field $V$ and $W$ of $N$ along $u$, i.e., 
\begin{equation}
\frac{\partial^2}{\partial s\partial t}\Bigg\vert_{(s,t)=(0,0)}E_2(u_{s,t}).
\end{equation}
We obtain:
\begin{thm} 
Let $u\,:(M,g)\rightarrow (N,h)$ be a statistical harmonic map, and 
$u_{s,t}$\,\,($\epsilon<s,t<\epsilon$), a smooth variation with two parameters $s$ and $t$ 
with $u_{0,0}=u$. 
 Then, it holds that 
 \begin{align}
 \frac{\partial^2}{\partial s\partial t} E_2(u_{s,t})
 \Bigg\vert_{(s,t)=(0,0)}\nonumber
 &=\int_M
 h\Big(W,-\sum_{i=1}^m\big({\widetilde{\nabla}}_{e_i} {\widetilde{\nabla}}_V u_{\ast}e_i
 - {\widetilde{\nabla}}_{ \nabla_{e_i} e_i }\big) V\\
 &\qquad\qquad\qquad  -\sum_{i=1}^mR^{N,h}(V,u_{\ast} e_i) u_{\ast} e_i
 \Big )\,d \mu_g\nonumber\nonumber\\
 &=\int_M h\Big(W, - \sum_{i=1}^m\big({\widetilde{\nabla}}_{e_i} {\widetilde{\nabla}}_{e_i} V-
 {\widetilde{\nabla}}_{ \nabla_{e_i} e_i } V\big)\nonumber\\
 &\qquad\qquad\qquad  -\sum_{i=1}^mR^{N,h}(V,u_{\ast} e_i) u_{\ast} e_i
 \Big )\,d \mu_g,
 \end{align}
 where $u_{s,t}$ is a smooth variation of $u$ with two parameters $s$ and $t$, such that 
 the mapping 
 $F\,;\,(-\epsilon, \epsilon)\times (-\epsilon,\epsilon)\times M\ni 
 (s,t,x)\mapsto u_{s,t}(x)\in N$ satisfies that 
 $$
 F_{\ast}\,\frac{\partial}{\partial\,t}\,(0,x)=V(x),\quad 
 F_{\ast}\,e_i(0,x)=u_{\ast}\,e_i(x),\quad i=1,\ldots,m.
 $$
 Here, $\{e_i\}_{i=1}^m$ is a locally defined orthonormal frame field 
 on $M$ with respect to the Riemannian metric $g$. 
\end{thm}
\vskip0.6cm\par
Here, let us recall the definition of the Jacobi operator which is the differential operator of the following form, 
\begin{defn}
\begin{align}\qquad\qquad J_uV&:=-\sum_{i=1}^m \Big({\widetilde{\nabla}}_{e_i} 
{\widetilde{\nabla}}_{e_i} V-
 {\widetilde{\nabla}}_{ \nabla_{e_i} e_i } \Big) \,V -
 \sum_{i=1}^mR^{N,h} ( V,u_{\ast} e_i) u_{\ast} e_i \nonumber\\
 &=\overline{\Delta}{}_uV-{\frak R}^u(V),\qquad\qquad\qquad\qquad\quad  V\in \Gamma(u^{-1}TN).
 \end{align}
 Here, for $V\in \Gamma(u^{-1}TN)$, we define the {\em rough Laplacian} 
 $\overline{\Delta}{}_u$ by
 \begin{equation}
 \overline{\Delta}{}_uV:=-\sum_{i=1}^m\Big({\widetilde{\nabla}}_{e_i}{\widetilde{\nabla}}_{e_i}
 -{\widetilde{\nabla}}_{\nabla_{e_i}e_i}\Big)V, 
 \end{equation}
 and the {\em curvature operator} $\frak R^u(V)$ by  
 \begin{equation}
 {\frak R}^u(V):= \sum_{i=1}^mR^{(N,h)}(V,u_{\ast}e_i)u_{\ast}e_i.
 \end{equation}
\end{defn}
\vskip0.3cm\par
{\it Proof of Theorem 3.1.} \quad 
Let us define $F(s,t,x):=u_{s,t}(x)\in M$ for 
$-\epsilon <s,t<\epsilon$ and $x\in M$ and small positive number $\epsilon>0$.  
Let us recall 
\begin{align*}
E_2(u_{s,t})&=\frac12\int_M\sum_{i=1}^mE_2(u_{s,t\,\ast}\,e_i, u_{s,t\,\ast}\,e_i)\, d\mu_g\\
&=\frac12 \int_M \sum_{i=1}^m h(F_{\ast}\,e_i, F_{\ast} \, e_i)\,d\mu_g, 
\end{align*}
Then, we have 
\begin{align*}
\frac{\partial}{\partial t}E_2(u_{s,t})&=\frac12\int_M\frac{\partial}{\partial t}\,
\sum_{i=1}^m h(F_{\ast}\,e_i.F_{\ast}\,e_i)\,d\mu_g\,\nonumber\\
&=\int_M\,\sum_{i=1}^m \, h(\widetilde{\nabla}_{\frac{\partial}{\partial t}}\,F_{\ast}\,e_i,\,F_{\ast}\,e_i)\,d\mu_g\nonumber\\
&=\int_M\bigg\{e_i\,\cdot\,h\bigg(
F_{\ast}\frac{\partial}{\partial\,t},F_{\ast}\,e_i\bigg)-
h\bigg(F_{\ast}
\frac{\partial}{\partial \,t},\,\widetilde{\nabla}_{e_i}F_{\ast}\,e_i\bigg)
\bigg\}\,d\mu_g\nonumber\\
&=\int_M\bigg\{{\rm div}(X_t)+\sum_{i=1}^m \bigg( h\big(F_{\ast}\,\frac{\partial}{\partial t},\,
F_{\ast}(\nabla_{e_i}e_i)\big)
-h\bigg(F_{\ast}\frac{\partial}{\partial t}, 
\widetilde{\nabla}_{e_i}F_{\ast}\, e_i
\bigg)\bigg)
\bigg\}\,d\mu_g\nonumber\\
&=-\int_M
h\bigg( F_{\ast}\frac{\partial}{\partial t},\,
\sum_{i=1}^m\bigg(\widetilde{\nabla}_{e_i}F_{\ast}e_i
-F_{\ast}\bigg(\nabla_{e_i}e_i
\bigg)
\bigg)\, d\mu_g.   \qquad\qquad  
\end{align*}
Since we have 
$$
F_{\ast}\,\frac{\partial}{\partial t}\,(0,x)=V(x),\,\,\,\, 
F_{\ast}\,e_i\,(0,x)=u_{\ast}e_i(x),\,\,\,\, 
F_{\ast}\,\nabla_{e_i}e_i\,(0,x)=u_{\ast}\nabla_{e_i}e_i(x),
$$
the right hand side coincides with 
\begin{align}
&-\int_Mh(V,\sum_{i=1}^m\big({\widetilde{\nabla}}_{e_i}\,u_{\ast}e_i
-u_{\ast}\nabla_{e_i}e_i
\big)\big)\,d\mu_g
=-\int_Mh(V,\tau(u))\,d\mu_g.
\end{align}
Here $\tau(u)\in \Gamma(u^{-1}TN)$ is the {\em tension field} of $u$ which is defined by 
\begin{equation}
\tau(u)(x):=\sum_{i=1}^m
\big(
{\widetilde{\nabla}}_{e_i}u_{\ast}e_i-u_{\ast}{\nabla}_{e_i}e_i
\big)(x).
\end{equation}
\vskip0.3cm\par
Next, we will derive the second variational formula. 
\par
Let $u \, : \, (M,g) \rightarrow (N,h)$ be a statistical harmonic mapping. 
Let $u_{s,t} \, : \, (M,g) \rightarrow (N,h)$ be a smooth variation with two parameters 
$s$ and $t$ of $u$, namely, let $F$ be a smooth mapping 
$$
F: \, (-\epsilon,\epsilon) \times (-\epsilon,\epsilon) \times M \ni (s,t,x) \mapsto u_{s,t}(x)\in N
$$
with 
$$
F(0,0,x)=u_{0,0}(x)=u(x),\qquad x \in M.
$$
Then, we will calculate the Hessian of the second variation of $E(u_{s,t})$. To do it, 
let $V$ and $W$ be two variation vector fields along $u$, which are defined by, for every 
$x \in M$, 
\begin{align}
V(x):= \frac{d}{ds} { \Big\vert }_{s=0} u_{s,0} (x) \in T_{u(x)} N,\quad 
W(x):= \frac{d}{dt} { \Big\vert }_{t=0} u_{0,t} (x) \in T_{u(x)} N, \nonumber
\end{align} 
respectively. 
\vskip0.3cm\par
We calculate 
\begin{equation}
H(E_2)_u(V,W)=\frac{\partial^2}{\partial s\partial t}{\Big\vert}_{(s,t)=(0,0)}E_2(u_{s,t}). 
\end{equation}
\vskip0.3cm\par
Indeed, we calculate 
\begin{align}
\frac{\partial^2}{\partial s\partial t}\,E_2(u_{s,t})&=-\int_M\frac{\partial}{\partial s}\,
h(F_{\ast}\frac{\partial}{\partial t}, \sum_{i=1}^m\big\{
\widetilde{\nabla}_{e_i} F_{\ast}e_i-F_{\ast}\nabla_{e_i}e_i\big\}\big)
d\mu_g\nonumber\\
&=-\int_Mh\bigg({\widetilde{\nabla}}_{\frac{\partial}{\partial s}}\,F{}_{\ast}
\frac{\partial}{\partial t},\,\sum_{i=1}^m
\big\{\widetilde{\nabla}_{e_i}F_{\ast}e_i
-F_{\ast}\nabla_{e_i}e_i\big\}
\bigg) \, d\mu_g\nonumber\\
&\quad -\int_Mh\bigg(F_{\ast}\frac{\partial}{\partial t},\,\sum_{i=1}^m
\widetilde{\nabla}_{\frac{\partial}{\partial s}}\big\{
\widetilde{\nabla}_{e_i}F_{\ast}e_i
-F_{\ast}\nabla_{e_i}e_i\big\}
\bigg)\,d\mu_g.
\end{align}
Here, the first term of the right hand side of $(3.9)$ must vanish. 
Because, in the integrand of the first term, 
$\sum_{i=1}^m
\big\{\widetilde{\nabla}_{e_i}F_{\ast}e_i
-F_{\ast}\nabla_{e_i}e_i\big\}$ coincides with 
$\tau(u)$ which must vanish at $(s,t)=(0,0)$ since we take for $u$ to be a 
statistical harmonic map. 
\vskip0.3cm\par
For the second term of $(3.9)$, we can see 
$$\widetilde{\nabla}_{\frac{\partial}{\partial s}}
\widetilde{\nabla}_{e_i}F_{\ast}e_i
=
\widetilde{\nabla}_{e_i}
\widetilde{\nabla}_{\frac{\partial}{\partial s}}
F_{\ast}e_i
+\widetilde{\nabla}_{
{\big[\frac{\partial}{\partial s}},
{e_i} 
\big]}
F_{\ast}e_i+R^{(N,h)}\bigg(F_{\ast}\frac{\partial}{\partial s}, 
F_{\ast}e_i\bigg)
F_{\ast}e_i. 
$$
But, it must hold that 
$\big[
\frac{\partial}{\partial s},
{e_i} 
\big]=0$. 
Furthermore, it holds that
\begin{align*}
\widetilde{\nabla}_{\frac{\partial}{\partial s}}
F_{\ast}
{\nabla}_{e_i}e_i
&=
\widetilde{\nabla}_{\nabla_{e_i}e_i}
F_{\ast}\frac{\partial}{\partial s}
+F_{\ast}
{\bigg[\frac{\partial}{\partial s}},
\nabla_{e_i} e_i
\bigg]\\
&=\widetilde{\nabla}_{\nabla_{e_i}e_i}
F_{\ast}\frac{\partial}{\partial s}\\
&=\widetilde{\nabla}_{\nabla_{e_i}e_i}V.
\end{align*}
Notice that 
at $(s,t)=(0,0)$, it holds that 
$$
F_{\ast}\frac{\partial}{\partial s}=V,\quad 
F_{\ast}\frac{\partial}{\partial t}=W,\quad
\mbox{and}\quad
F_{\ast}e_i=u_{\ast}e_i. 
$$
\vskip0.3cm\par
Therefore, we obtain
$$\left\{
\begin{aligned}
\widetilde{\nabla}_V\widetilde{\nabla}_{e_i}u_{\ast}e_i
&=\widetilde{\nabla}_{e_i}\widetilde{\nabla}_Vu_{\ast}e_i+R^{(N,h)}(V,u_{\ast}e_i)u_{\ast}e_i,\\
\widetilde{\nabla}_V\nabla_{e_i}e_i&={\widetilde{\nabla}}_{\nabla_{e_i}e_i}V.
\end{aligned}
\right.
$$
Therefore, we obtain 
\begin{align}
\frac{\partial^2}{\partial s\partial t}E_2(u_{s,t})
\bigg\vert_{(s,t)=(0,0)}
&=
-\int_Mh\bigg(W,\sum_{i=1}^m\bigg\{\widetilde{\nabla}_{e_i}{\widetilde{\nabla}}_Vu_{\ast}e_i
+R^{(N,h)}(V,u_{\ast}e_i)u_{\ast}e_i-\widetilde{\nabla}_{{\nabla_{e_i}}e_i}V\bigg\}\bigg)\,d\mu_g\nonumber\\
&=-\int_Mh\bigg(W,\sum_{i=1}^m\bigg\{
{\widetilde{\nabla}}_{e_i}{{\widetilde{\nabla}}_{e_i}}V
-{\widetilde{\nabla}}_{\nabla_{e_i}e_i}V
+R^{(N,h)}(V,u_{\ast}e_i)u_{\ast}e_i
\bigg\}\bigg)\,d\mu_g.
\end{align}
In (3.10), we only have to see 
\begin{align}
\widetilde{\nabla}_Vu_{\ast}e_i=\widetilde{\nabla}_{e_i}V,
\end{align}
because it holds that the torsion tensor field 
$T(X,Y)$ holds that 
\begin{align*}
0=T(V,e_i)
&=\widetilde{\nabla}_V u_{\ast}e_i-{\widetilde{\nabla}}_{e_i}V-{[V,e_i]}\\
&=\widetilde{\nabla}_V u_{\ast}e_i-{\widetilde{\nabla}}_{e_i}V. 
\end{align*}
\vskip0.3cm\par
Therefore, we obtain 
\begin{equation}
\frac{\partial^2}{\partial s\partial t}E_2(u_{s,t})\bigg\vert_{(s,t)=(0,0)}=\int_Mh(W, J_u(V))\,d\mu_g, 
\end{equation}
where the operator $J_u$, called the {\em Jacobi operator}, is given by 
\begin{equation}
J_u(V)=-\sum_{i=1}^m\bigg(
{\widetilde{\nabla}}_{e_i}
{\widetilde{\nabla}}_{e_i}V-{\widetilde{\nabla}}_{\nabla_{e_i}e_i}V
\bigg)
-\sum_{i=1}^mR^{(N,h)}(V,u_{\ast}e_i)u_{\ast}e_i.
\end{equation}
Thus, we can restate Theorem 3.1 as follows: 
\vskip0.3cm\par
\begin{thm} 
Let $u\,:(M,g)\rightarrow (N,h)$ be a statistical harmonic map, and 
$u_{s,t}$\,\,($\epsilon<s,t<\epsilon$), a smooth variation with two parameters $s$ and $t$ 
with $u_{0,0}=u$. 
 Then, it holds that 
 \begin{align}
 \frac{\partial^2}{\partial s\partial t} E_2(u_{s,t})
 \Bigg\vert_{(s,t)=(0,0)}
 &=\int_M
 h\Big(W,J_u(V)
 \Big )\,d \mu_g,
 \end{align}
 where $u_{s,t}$ is a smooth variation of $u$ with two parameters $s$ and $t$, such that 
 the mapping 
 $F\,;\,(-\epsilon, \epsilon)\times (-\epsilon,\epsilon)\times M\ni 
 (s,t,x)\mapsto u_{s,t}(x)\in N$ satisfies that 
 $$
 F_{\ast}\,\frac{\partial}{\partial\,t}\,(0,x)=V(x),\quad 
 F_{\ast}\,e_i(0,x)=u_{\ast}\,e_i(x),\quad i=1,\ldots,m.
 $$
 and 
 $J_u$ is the Jacobi operator defined by $(3.13)$. 
 Here, $\{e_i\}_{i=1}^m$ is a locally defined orthonormal frame field 
 on $M$ with respect to the Riemannian metric $g$. 
\end{thm}
\vskip0.6cm\par
\section{The index, nullity and the weak stability of a statistical harmonic mapping.}
In this section, 
we want to discuss the index, nullity and the weak stability for a statistical harmonic mapping 
analogously  as the usual harmonic mapping (cf. \cite{29}). 
 \vskip0.3cm\par
Let us denote the spectrum of the operator 
$J_u$ by 
\begin{align}
\lambda_1(u)\leqq\lambda_2(u)\leqq \cdots\leqq\lambda_i(u)\leqq \cdots.
\end{align}
counted the eigenvalues with their multiplicities.  
\par
Let us denote by $V_{\lambda}$, the {\it eigenspace} of the operator $J_u$ with the {\it eigenvalue} $\lambda$ 
when
\begin{align}
V_{\lambda}&:=\left\{
V\in \Gamma(u^{-1}TN):\,J_uV=\lambda\,V
\right\}\\
&\not=\{0\}.\nonumber
\end{align}
\par 
The {\it index}  of the eigenvalue $\lambda$, denoted by {\rm index}$({\lambda})$, is defined by 
\begin{align}
{\rm sup}\{{\rm dim}(F):\,& F\subset\Gamma(u^{-1}TN, 
{\rm\,  a \, subspace} \, \nonumber \\ 
\qquad \qquad 
{\rm on \,\, which \,\,} & H(E_2)_u {\rm \, is \, \, negative \, \, definite}
\}, 
\end{align}
\par 
The {\it nullity} of the eigenvalue $\lambda$, denoted by nullity$(\lambda)$, is defined by 
\begin{align}
{\rm dim} \{
V\in \Gamma(u^{-1}TN):\,\,H(E_2)_{u}(V,W)=\{0\}\,\,\, \forall\,\,W \in \Gamma(u^{-1}TN)
\}. 
\end{align}
\par
\par 
\begin{defn}
A statistical harmonic mapping $u:\,(M,g)\rightarrow (N,h)$ is said to be 
\it{weakly stable} 
{\rm if} index$(u)=0$. 
{\rm Otherwise, it is said to be} \it{unstable}. 
\end{defn}
\par
\begin{thm}
Assume that a statistical Riemannian manifold 
$(N,h,\nabla)$ is non-positive curvature, i.e., $h(R^{(N,h)}(U,V)V,U)\leqq 0$ 
for all $U$ and $V$ in $\Gamma(u^{-1}TN)$. 
Then, a statistical harmonic mapping $u:\,(M,g,\nabla^M)\rightarrow (N,h,\nabla^N)$ is weakly stable, i.e., it holds that 
. \begin{align}
\int_Mh(J_u(V),V)\,d\mu_g\geqq 0\qquad\qquad (\forall\,V\in \Gamma(\mu^{-1}TN)).
\end{align}
\end{thm}
{\it Proof.} \rm 
Since a statistical Riemannian manifold 
$(N,h,\nabla)$ is non-positive, 
it holds that,  
for all $i=1,\cdots, m$, and $V\in \Gamma(\mu^{-1}TN)$, 
\begin{align}
\int_M h(R^{(N,h)}(V,u_{\ast}
e_i)u_{\ast}e_i,V)d\mu_g\leqq 0. 
\end{align}
Furthermore, 
the rough Laplacian $\overline{\Delta}$ satisfies that 
\begin{align*}
(
{\overline{\Delta}}_u
V,V
)
&=-\int_M
\bigg(
\sum_{i=1}^n
\bigg\{
{\widetilde{\nabla}}_{e_i}
{\widetilde{\nabla}}_{e_i}V
-
{\widetilde{\nabla}}_
{
{{\nabla}_{e_i}}
e_i}V
\bigg\},V
\bigg)\,
d\mu_g\\
&=\sum_{i=1}^m
\int_M
\bigg(
{\widetilde{\nabla}}_{e_i}V,
{\widetilde{\nabla}}_{e_i}V
\bigg)
d\mu_g\\
&\geqq 0,
\end{align*}
and 
\begin{align*}
(
{\overline{\Delta}}_u
V,V
)=0
\qquad &{\Longleftrightarrow\qquad 
\widetilde{\nabla}}_{e_i}V=0\qquad (\forall \,\, i=1,\cdots,m)\\
&\Longleftrightarrow \qquad {\widetilde{\nabla}}_XV=0\qquad 
(\forall \,\, X \in {\frak X}(M))\\
&{\Longleftrightarrow}\quad \,\,
\tau(u)=0.
\end{align*}
Thus, $u:(M,g)\rightarrow (N,h)$,\,{\rm statistically harmonic}. 
Therefore, 
$u:\,(M,g,\nabla^M)\rightarrow (N,h,\nabla^N)$ is weakly stable. 
\qed
\vskip0.6cm\par
Furthermore, the assumption of this theorem can be weakened into the following way: 
\begin{defn}
Let $u:\,(M,g,\nabla)\rightarrow (N,h,\nabla)$ be a smooth map of a compact statistical manifold into another statistical manifold. 
Then a Riemannian manifold $(N,h)$ is 
 {\em non-positive in the sense of integration} 
if it holds that for every $U$ and $V$ belonging to $\Gamma(\mu^{-1}TN)$, 
\begin{equation}
\int_Mh(R^{(N,h)}(U,V)V,U)\,d\mu_g\leqq 0.
\end{equation}
\end{defn}
\vskip0.6cm\par 
Thus, we obtain the following theorem. 
\begin{thm}
Assume that a Riemannian manifold $(N,h)$ is 
 {\em non-positive in the sense of integration}.  
Then, every smooth harmonic map $u:\,(M,g,\nabla)\rightarrow (N,h,\nabla)$ of a compact statistical manifold into another statistical manifold 
is weakly stable. 
\end{thm}
The proof can be given in the same way as Theorem 4.2, so it is omitted. \qed 
\vskip0.6cm\par
 \section{Examples}
 In this section, we will give several examples of statistical harmonic mappings. 
 \begin{ex}(the constant maps) 
 Let $(M,g,\nabla)$ be a compact Riemannian manifold. 
 Then every constant mapping 
 $\phi:\,M\rightarrow N$ is weakly stable statistical harmonic, and  
\vskip0.3cm \par\qquad
 index$(\phi)=0$,\qquad   and \qquad  nullity$(\phi)=\dim(N)$.  
 \vskip0.3cm\par
 \end{ex}
 \begin{ex}(the identity maps) Let $(M,g)$ be a compact Einstein Riemannian manifold, i.e., 
 Ric${}_M=c\,g$ with some constant $c$. Let id$:\,M\rightarrow M$ be the identity mapping. 
 Then, it holds that 
 \vskip0.3cm\par\qquad 
 index$($id$)=0\qquad\qquad \Longleftrightarrow \qquad\qquad 
 \lambda_1(g)\geqq 2c.$ 
 \vskip0.3cm\par
 \end{ex}
 \begin{ex}
 Let $(S^n,{\rm id})$ be the unit sphere with constant sectional curvature $1$ $(n\geqq 3)$. 
Let $\phi:\,(M^m,g)\rightarrow (S^n,{\rm id})$ be an arbitrary non-constant harmonic 
mapping. Then, \quad 
index$(\phi)>0$. 
 \vskip0.3cm\par
 \end{ex}
 \begin{ex}
 Let $(N,h)$ be a compact Riemannian manifold with non-positive sectional curvature, 
 and $\phi:\,(M,g)\rightarrow (N,h)$, an arbitrary harmonic mapping from an arbitrary compact Riemannian manifold into this Riemannian manifold $(N,h)$. 
 Then, \quad 
  index$(\phi)=0. $
 \end{ex}
 \vskip0.3cm\par

\end{document}